\documentclass[a4paper,12pt]{amsart}
\usepackage{amssymb}
\usepackage{ifthen}
\usepackage[dvips]{graphicx}
\usepackage{cite}
\nonstopmode \numberwithin{equation}{section}
\usepackage{amssymb}
\usepackage{ifthen}
\usepackage{graphicx}
\usepackage{amsmath}
\usepackage[T1]{fontenc} 
\usepackage[utf8]{inputenc}
\usepackage[usenames,dvipsnames]{color}
\usepackage{color}
\usepackage[english]{babel}
\usepackage{fancyhdr}
\usepackage{fancybox}
\usepackage{tikz}

\nonstopmode \numberwithin{equation}{section}
\theoremstyle{plain}

\newtheorem{prop}{Proposition}

\newtheorem{conj}{Conjecture}

\theoremstyle{definition}
\newtheorem{defn}{Definition}[section]
\newtheorem{example}{Example}[section]
\newtheorem{thm}{Theorem}[section]
\newtheorem{cor}{Corollary}[section]
\newtheorem{lem}{Lemma}[section]
\newtheorem{prob}{Problem}
\newtheorem{rem}{Remark}[section]
\newtheorem{ques}{Question}[section]

\newcounter{minutes}
\divide\time by 60
\newcounter{hours}
\multiply\time by 60
\addtocounter{minutes}{-\time}

\newcounter {own}
\def\theown {\thesection       .\arabic{own}}

\newenvironment{pf}[1][]{%
 \vskip 3mm
 \noindent
 \ifthenelse{\equal{#1}{}}%
  {{\slshape Proof. }}%
  {{\slshape #1.} }%
 }%
{\qed\bigskip}

\newcounter{alphabet}

\def\be{\begin{equation}}
\def\ee{\end{equation}}

\newcommand{\bee}{\begin{enumerate}}
\newcommand{\eee}{\end{enumerate}}

\newcommand{\blem}{\begin{lem}}
\newcommand{\elem}{\end{lem}}
\newcommand{\bthm}{\begin{thm}}
\newcommand{\ethm}{\end{thm}}
\newcommand{\bcor}{\begin{cor}}
\newcommand{\ecor}{\end{cor}}
\newcommand{\beg}{\begin{examp}}
\newcommand{\eeg}{\end{examp}}
\newcommand{\begs}{\begin{examples}}
\newcommand{\eegs}{\end{examples}}

\newcommand{\bdefn}{\begin{defn}}
\newcommand{\edefn}{\end{defn}}

\newcommand{\bprob}{\begin{prob}}
\newcommand{\eprob}{\end{prob}}
\newcommand{\bei}{\begin{itemize}}
\newcommand{\eei}{\end{itemize}}

\newcommand{\bcon}{\begin{conj}}
\newcommand{\econ}{\end{conj}}
\newcommand{\bcons}{\begin{conjs}}
\newcommand{\econs}{\end{conjs}}
\newcommand{\bprop}{\begin{prop}}
\newcommand{\eprop}{\end{prop}}
\newcommand{\br}{\begin{rem}}
\newcommand{\er}{\end{rem}}
\newcommand{\brs}{\begin{rems}}
\newcommand{\ers}{\end{rems}}
\newcommand{\bo}{\begin{obser}}
\newcommand{\eo}{\end{obser}}
\newcommand{\bos}{\begin{obsers}}
\newcommand{\eos}{\end{obsers}}
\newcommand{\bpf}{\begin{pf}}
\newcommand{\epf}{\end{pf}}
\newcommand{\ba}{\begin{array}}
\newcommand{\ea}{\end{array}}
\newcommand{\beq}{\begin{eqnarray}}
\newcommand{\beqq}{\begin{eqnarray*}}
\newcommand{\eeq}{\end{eqnarray}}
\newcommand{\eeqq}{\end{eqnarray*}}

\begin{document}

\title{Entire Solutions for quadratic trinomial-type partial differential-difference equations in $ \mathbb{C}^n $}

\author{Sanju Mandal}
\address{Sanju Mandal,
	Department of Mathematics,
	Jadavpur University,
	Kolkata-700032, West Bengal, India.}
\email{sanjum.math.rs@jadavpuruniversity.in}

\author{Molla Basir Ahamed}
\address{Molla Basir Ahamed,
	Department of Mathematics,
	Jadavpur University,
	Kolkata-700032, West Bengal, India.}
\email{mbahamed.math@jadavpuruniversitry.in}

\subjclass[{AMS} Subject Classification:]{Primary 39A45, 30D35, 35M30, 32W50}

\def\thefootnote{}
\footnotetext{ {\tiny File:~\jobname.tex,
printed: \number\year-\number\month-\number\day,
          \thehours.\ifnum\theminutes<10{0}\fi\theminutes }
} \makeatletter\def\thefootnote{\@arabic\c@footnote}\makeatother

\begin{abstract}
 In this paper, utilizing Nevanlinna theory, we study existence and forms of the entire solutions $ f $ of the quadratic trinomial-type partial differential-difference equations in $ \mathbb{C}^n $
 \begin{align*}
 	a\left(\alpha\dfrac{\partial f(z)}{\partial z_i} + \beta\dfrac{\partial f(z)}{\partial z_j}\right)^2 + 2 \omega \left(\alpha\dfrac{\partial f(z)}{\partial z_i} + \beta\dfrac{\partial f(z)}{\partial z_j}\right) f(z + c)  + b f(z + c)^2 = e^{g(z)}
 \end{align*}
and
 \begin{align*}
 	a\left(\alpha\dfrac{\partial f(z)}{\partial z_i} + \beta\dfrac{\partial f(z)}{\partial z_j}\right)^2 & + 2 \omega \left(\alpha\dfrac{\partial f(z)}{\partial z_i} + \beta\dfrac{\partial f(z)}{\partial z_j}\right) \Delta_cf(z) + b [\Delta_cf(z)]^2 = e^{g(z)},
 \end{align*}
where $ a, \omega, b\in\mathbb{C} $, $ g $ is a polynomial in $ \mathbb{C}^n $ and $ \Delta_cf(z)=f(z+c)-f(z) $. The main results of the paper improve several existence results in $ \mathbb{C}^n $ for integer $ n\geq 2 $ and $ 1\leq i<j\leq n $  and their corollaries of the paper are an extension of the results of Xu \emph{et al. } [\textit{Rocky Mountain J. Math.} \textbf{52}(6) (2022), 2169–2187] for trinomial equation with arbitrary coefficient in $ \mathbb{C}^2 $. Moreover, examples are exhibited to validate the conclusion of the main results.
\end{abstract}

\maketitle
\pagestyle{myheadings}
\markboth{Molla Basir Ahamed and Sanju Mandal}{Entire solutions for trinomial-type equations in $ \mathbb{C}^n $}

\section{Introduction}
In this paper, we consider meromorphic solutions of certain functional equations in $ \mathbb{C}^n $ related to Fermat varieties. Among the most basic functional equations are the circle functional equation $ f^2 + g^2 = 1 $, and the Fermat cubic $ f^3 + g^3 = 1 $. Generalizations of these power equations are called Fermat-type functional equations, which are associated with diagonal varieties, and have been the subject of interest in global complex analysis in connection with the extensions of Picard-type theorems and results on hyperbolic sub-manifolds of projective space (see for example \cite{Green-AJM-1975,Kiernan-PAMS-1969, Yang & 1970}). Due to the development of the difference analogue lemma of logarithmic derivative lemma, in recent year an increasing amount of interests has been grown up for several properties of entire and meromorphic solutions of several difference functional equations both in one and several complex variables. Since non-constant polynomials in $ \mathbb{C}^n $ (for $ n\geq 2 $) may be periodic, the nature of solutions of Fermat-type equations in $ \mathbb{C}^n $ is completely different from that in $ \mathbb{C} $. This is one of the reason why we consider Fermat-type functional equations in several complex variables in our study. \vspace{1.2mm}

The study of Fermat-type functional equation has been an interesting subject in the field of complex analysis in connection with extensions of Nevanlinna's theory. For extensive research on the Fermat-type functional equations, we refer to the articles \cite{Cao-MN-2013,Cao-Xu-ADM-2020,Chen-JMAA-2022, XU-CAO-MJM-2018,Xu-Liu-Li-JMAA-2020,Xu-RMJM-2021,Xu-ZHANG-ZHENG-RMJM-2022} and references therein. We will assume that the reader is familiar with basic elements of the Nevanlinna's theory of meromorphic function $ f $ in one or several complex variables (see e.g., \cite{Hayman-CP-1964,Yang-SSP-1993, Hu-Li-Yang-2003,Ye-MZ-1996}), such as the characteristic function $ T(r,f) $, the counting function $ N(r,f) $ for poles of $ f $, reduce counting function $ \overline{N}(r,f) $ of $ f $, proximation function $ m(r,f) $ in the value distribution theory, also known as Nevanlinna theory. We denote by $ S(r, f) $, any function satisfying $ S(r, f) = \circ\{T(r,f)\} $  as $ r\rightarrow\infty $, possibly outside a set of finite measure. In addition, we use the notation $ \rho(f) $ to denote the order of growth of the meromorphic function $ f $ in $ \mathbb{C}^n $, and defined by
\begin{align*}
	\rho(f) =\limsup_{r\rightarrow\infty} \dfrac{\log T(r,f)}{\log r}.
\end{align*}

It has always been a well-known and interesting problem to investigate the existence and form of solutions to Fermat-type functional equations of the form
\begin{align}\label{eq-1.1}
	f^n(z)+g^n(z)=1
\end{align}
regard as the Fermat diophantine equation $ x^n+y^n=1 $ over functional fields, where $ n\geq 2 $ is an integer. The classical results on meromorphic solutions in $ \mathbb{C} $ of \eqref{eq-1.1} have been studied and forms of the solutions are obtained (see e.g. \cite{Baker-PAMS-1966,Gross-BAMS-1966,Montel & Paris & 1927}). It is understood that \eqref{eq-1.1} does not admit transcendental meromorphic (resp. entire) solutions when $ n\geq 4 $ (resp. $ n\geq 3 $). If $ n=3 $, then equation \eqref{eq-1.1} admits meromorphic solutions $ f ={(3 +\sqrt{3}\wp^{\prime}(\beta))}/{6\wp(\beta)} $ and $ g=\eta{(3 -\sqrt{3}\wp^{\prime}(\beta))}/{6\wp(\beta)} $, for some non-constant entire function $ \beta $, where $ \eta^3 =1 $ and $ \wp $ denotes the Weierstrass $ \wp $-function satisfying $ (\wp^{\prime})^2 \equiv 4\wp^3 -1 $ after appropriately choosing its periods. For $ n=2 $, \eqref{eq-1.1} has nontrivial (non-constant) entire solutions $ f(z)=\cos(\psi(z)) $ and $ g(z)=\sin(\psi(z)) $, where $ \psi $ is an entire function. For the study of meromorphic solutions to \eqref{eq-1.1}
in $ \mathbb{C}^n $ and applications to complex partial differential
equations, we refer the reader to (see e.g \cite{Li-FM-2005,Li-IJM-2014,Li-NMJ-2005,Li-AM-2007}). \vspace{1.2mm}

This article mainly concerns the global analytic or meromorphic solutions for trinomial quadratic partial differential-difference equations (in short, $ PDDE $) with arbitrary coefficients of the form $ aF^2+2\omega FG+bG^2=e^g $, where $ a, b, \omega $ are complex constants and $ g $ is a polynomial in $ \mathbb{C}^n $. In particular, for $ a=1=b $ and $ g(z)=2k\pi i $, $ k $ being an integer, there are number of results in $ \mathbb{C} $ and $ \mathbb{C}^2 $. In fact, what could be the characterization of solutions of the trinomial in $ \mathbb{C}^n $ is not explored yet and need to study. In general, one cannot expect the existence of analytic solutions, and even when global analytic or entire solutions exist, it is difficult to find such solutions in closed form in $ \mathbb{C}^n $. The finite order solutions to the Fermat-type binomial and trinomial equations in $ \mathbb{C} $ over some commonly studied function fields have been investigated by many authors, and there is an extensive literature on these equations and generalizations as well as connections to other problems (see e.g., \cite{Cao-MN-2013,Gross-BAMS-1966,Gundersen-CMFT-2017,Montel & Paris & 1927,Tang & Liao & 2007,Yang & 1970,Yang & Li & 2004,Gundersen-Heyman-2004-BLMS}). Furthermore, it appears that the solutions of the system of Fermat-type binomial or trinomial equations in $ \mathbb{C}^2 $ has been recently studied in \cite{Xu-RMJM-2021, Xu-Liu-Li-JMAA-2020}. However, no study has so far been done on the solutions of quadratic trinomial functional equations in $ \mathbb{C}^n$. In this paper, our main aim is to describe transcendental solutions for quadratic trinomial $ PDDEs $ in $ \mathbb{C}^n $. \vspace{1.2mm}

Liu \emph{et al.} \cite{Liu-Cao-Cao-AM-2012} have investigated the Fermat-type difference equation $ f^2 (z) + f^2 (z +c) = 1 $ in $ \mathbb{C} $ and obtained the finite order transcendental entire solutions satisfy $ f(z)=\sin(Az+B) $, where $ B $ is a constant and $ A = ((4k+1)\pi)/(2c) $, where $ k $ is an integer. Later, Han and Lü \cite{Han-JCMA-Lu-2019} established the solution to the more general complex difference equation  $ f^n(z) +g^n(z) = e^{\alpha z+\beta} $. Moreover, Liu \emph{et al.} \cite{Liu-Cao-Cao-AM-2012} showed that the existence of solutions for the complex differential-difference equations $ f^{\prime} (z)^2 + f(z +c)^2 = 1 $ and $ f^{\prime} (z)^2 + [f(z +c) -f(z)]^2 = 1 $ in $ \mathbb{C} $. \vspace{1.2mm}

As is known to all, partial differential equations (PDEs) are occurring in various areas of applied mathematics, such as fluid mechanics, nonlinear acoustics, gas dynamics, and traffic flow (see \cite{Courant-Hilbert-I-1962, Garabedian-W-1964}). In general, it is difficult to find entire and meromorphic solutions for a nonlinear PDE. By employing Nevanlinna theory and the method of complex analysis, there were a number of literature focusing on the solutions of some PDEs and theirs many variants, readers can refer to \cite{Chen-JMAA-2022,Cao-Xu-ADM-2020,Khavinson & Am. Math. Mon & 1995,Li-TAMS-2004,Li-AM-2007, Lu-CRM-2020,Gundersen-PEMS-2020, Saleeby-AM-2013,Xu-AMP-2022}.\vspace{1.2mm}

The solutions of Fermat-type $ PDEs $ were investigated by \cite{Li-IJM-2004,Saleeby-A-1999}. Most noticeably, in $ 1995 $, Khavinson \cite{Khavinson & Am. Math. Mon & 1995} derived that any entire solution of the partial differential equation in $ \mathbb{C}^2 $,
\begin{align*}
	\left(\frac{\partial u}{\partial z_1}\right)^2+\left(\frac{\partial u}{\partial z_2}\right)^2=1
\end{align*}
is necessarily linear, i.e., $ u(z_1, z_2)=az_1+bz_2+c $, where $ a,b,c\in \mathbb{C} $, and $ a^2+b^2=1 $. This $ PDE $ in the real variable case occurs in the study of characteristic surfaces and wave propagation theory, and it is the two-dimensional eiconal equation, one of the main equations of geometric optics (see \cite{Courant-Hilbert-I-1962}). Furthermore, Li \cite{Li-NMJ-2005,Li-AM-2007} have continued the research and discussed solutions of a series of $ PDEs $ with more general forms including $ \left(\frac{\partial f}{\partial z_1}\right)^2 + \left(\frac{\partial f}{\partial z_2}\right)^2 =e^g $, $ \left(\frac{\partial f}{\partial z_1}\right)^2 + \left(\frac{\partial f}{\partial z_2}\right)^2 = p $, etc., where $ g, p $ are polynomials in $ \mathbb{C}^2 $. Recently, Xu \emph{et al.} \cite{Xu-ZHANG-ZHENG-RMJM-2022} established the solution of the $ PDDEs $
\begin{align}\label{eq-1.2}
 	\left(\alpha\frac{\partial f(z)}{\partial z_1} + \beta\frac{\partial f(z)}{\partial z_2}\right)^2 +  f(z + c)^2 = e^{g(z)}
\end{align}
and 
\begin{align}\label{eq-1.3}
	\left(\alpha\frac{\partial f(z)}{\partial z_1} + \beta\frac{\partial f(z)}{\partial z_2}\right)^2 +  [f(z + c)-f(z)]^2 = e^{g(z)}
\end{align}
in $ \mathbb{C}^2 $, and they obtained the form of the solution in $ \mathbb{C}^2 $.\vspace{1.2mm}

Inspired by the above results a question can be raised naturally:
\begin{ques}
What can be said about the form of solutions in $ \mathbb{C}^n $, if we extend the binomial equations \eqref{eq-1.2} and \eqref{eq-1.3} in \cite[Theorem 2.1, Theorem 2.2] {Xu-ZHANG-ZHENG-RMJM-2022} to trinomial equation with arbitrary coefficient.
\end{ques}

Motivated by the above question, our purpose of this article is to exploring the finite order transcendental entire solutions of the quadratic trinomial
partial differential equations. To find precise solutions of trinomial quadratic functional equations we use with certain techniques. More precisely,  Saleeby \cite{Saleeby-AM-2013} initiates this type of study considering the quadratic trinomial equations of the form 
$ f^2+2\alpha fg+g^2=1, $
where $ \alpha\in\mathbb{C}\setminus\{-1, 1\} $, which is associated with the partial differential equations 
\begin{align}\label{eq-1.4}
	u_{x}^2 + 2\alpha u_xu_y + u_y^2 = 1,
\end{align} 
where $ (x, y)\in\mathbb{C}^2 $ and showed that the entire and meromorphic solutions of \eqref{eq-1.4} have the form $ u(x,y)=ax+by+c $, where $ a^2 + 2\alpha ab + b^2=1 $.\vspace{1.2mm} 

The main tools are used in this paper are the Nevanlinna theory and the characteristic
equations for quasi-linear $ PDEs $ and linear $ PDEs $. The paper is organized as follows. Our main results about the existence and the forms of entire solutions and their corollaries with examples will be exhibited in Section 2. The proofs of the main results will be given in Section 3. 

\section{Main results}
Motivated by method of proof of results in \cite{Xu-ZHANG-ZHENG-RMJM-2022}, we explore the finite order transcendental entire solutions of quadratic trinomial partial differential equations in $ \mathbb{C}^n $. Henceforth, throughout this paper, we assume that $ z+c=(z_1 + c_1, \ldots, z_n + c_n) $, for any $ z=(z_1,\ldots,z_n) $ and $ c=(c_1,\ldots,c_n) $ are in $ \mathbb{C}^n $. To serve the purpose, we define  $ \omega_1:=-\frac{\omega}{\sqrt{ab}} \pm \frac{\sqrt{\omega^2 -ab}}{\sqrt{ab}} $ and $ \omega_2: =-\frac{\omega}{\sqrt{ab}} \mp \frac{\sqrt{\omega^2 -ab}}{\sqrt{ab}} $. Let $ g(z)=\sum_{|I|=0}^{p} a_{\alpha_1,\ldots,\alpha_n} z^{\alpha_1}_1\cdots z^{\alpha_n}_n $ be polynomial in $ \mathbb{C}^n $, where $ I=(\alpha_1,\ldots,\alpha_n) $ be two multi-index with $ |I|= \sum_{j=0}^{n}\alpha_j $ and $ \alpha_j $ are non-negative integers.

\vspace{1.2mm}

With the help of a transformation in trinomial $ PDDEs $, we obtain the following result concerning existence and forms of the solutions $ f $ of
\begin{align}\label{eq-2.1}
	a\left(\alpha\dfrac{\partial f(z)}{\partial z_i} + \beta\dfrac{\partial f(z)}{\partial z_j}\right)^2 + 2 \omega \left(\alpha\dfrac{\partial f(z)}{\partial z_i} + \beta\dfrac{\partial f(z)}{\partial z_j}\right) f(z + c)  + b f(z + c)^2 = e^{g(z)}.
\end{align} 
\begin{thm}\label{th-2.1}
Let $ c\in\mathbb{C}^n\setminus\{0\} $, $ a,b,\alpha\neq 0 $, and $ \omega^2\neq ab $. For $ 1\leq i< j\leq n $, if the $ PDDE $ \eqref{eq-2.1} in $ \mathbb{C}^n $ admits a transcendental entire solution of finite order, then $ g $ must be a polynomial of the form $ g(z)= L(z) + H(s) + B_1 $, where $ L(z)= a_1 z_1 +\cdots + a_n z_n $ and $ H(s) $ is a polynomial in $ s:= d_1 z_1 +\cdots + d_n z_n $ in $ \mathbb{C}^n $ with $ d_1 c_1 +\cdots + d_n c_n = 0 $ with $ H(z+c)= H(z) $, $ a_1,\ldots,a_n, B_1\in \mathbb{C} $. Furthermore, $ f $ must assume one of the following forms:
\begin{enumerate}
	\item [(i)] 
	\begin{align*}
		f(z)= \pm \sqrt{\dfrac{(b\omega^2_2 + a)\omega^2_1 - 2(\omega_1\omega_2 b + a)\omega_1\omega_2 + (b\omega^2_1 + a)\omega^2_2} {ab\omega_1\omega_2(\omega_2 -\omega_1)^2}} e^{\frac{g(z-c)}{2}},
	\end{align*}
    where $ g(z)= \psi(z_j - \frac{\beta}{\alpha}z_i) $ and $ \psi $ is a polynomial in $ \mathbb{C}^n $;
    \item [(ii)] 
    \begin{align*}
    	f(z) = \dfrac{\xi^2 -1}{\xi\sqrt{b}(\omega_2-\omega_1)} e^{\frac{L(z)+ H(s) + B}{2}}, \; \mbox{where}\; B\in\mathbb{C}
    \end{align*}
    with
    \begin{align*}
    	\dfrac{\sqrt{a}(\xi^2 -1)}{2\sqrt{b}(\omega_2\xi^2 -\omega_1)}(\alpha a_{i} + \beta a_{j}) = e^{\frac{a_1 c_1 +\cdots + a_n c_n}{2}};
    \end{align*}
    \item [(iii)] 
    \begin{align*}
    	f(z) = \dfrac{e^{L_1(z) +H_1(s) + E_1 -L_1(c)} - e^{L_2(z) +H_2(s) + E_2 -L_2(c)}}{\sqrt{b}(\omega_2 -\omega_1)},
    \end{align*}
    where $ L_{l}(z)= a_{l1} z_1 +\cdots + a_{ln} z_n $ for $ l=1, 2 $ with $ E_1, E_2\in\mathbb{C} $ and $ H_{l}(s) $ (for $ l=1, 2 $) are polynomials in $ s:= d_1 z_1 +\cdots + d_n z_n $ in $ \mathbb{C}^n $ with $ d_1 c_1 +\cdots + d_n c_n = 0 $ and $ H_{l}(z+c)= H_{l}(z) $ for $ l=1, 2 $, such that
    \begin{align*}
    	L_1(z) +H_1(s)\neq L_2(z) +H_2(s), \;\; g(z)= L_1(z) + L_2(z) + H_1(s) + H_2(s) + E_1 + E_2,
    \end{align*}
    and
    \begin{align*}
    	\dfrac{\sqrt{a}}{\omega_2\sqrt{b}}[\alpha a_{1i} +\beta a_{1j}] e^{-L_1(c)}\equiv 1 \;\;\; \mbox{and} \;\;\; \dfrac{\sqrt{a}}{\omega_1\sqrt{b}}[\alpha a_{2i} +\beta a_{2j}] e^{-L_2(c)}\equiv 1.
    \end{align*}
\end{enumerate}	
\end{thm}

\begin{rem}
Theorem \ref{th-2.1} is an extension of the result of Xu et al. \cite[Theorem 2.1]{Xu-ZHANG-ZHENG-RMJM-2022} and Xu et al. \cite[Theorem 1.2]{XU-CAO-MJM-2018} in $ \mathbb{C}^n $.
\end{rem}

The following result is an immediate corollary of Theorem \ref{th-2.1} for the solution to the trinomial partial differential-difference equations in $ \mathbb{C}^2 $ and this result can be considered as a trinomial version with arbitrary coefficients of that binomial equation \eqref{eq-1.2} in \cite{Xu-ZHANG-ZHENG-RMJM-2022}

\begin{cor}\label{cor-2.1}
Let $ c\in\mathbb{C}^2\setminus\{0\} $, $ a,b,\alpha\neq 0 $, and $ \omega^2\neq ab $. If the partial differential-difference equation	
\begin{align*}
	a\left(\alpha\dfrac{\partial f(z)}{\partial z_1} + \beta\dfrac{\partial f(z)}{\partial z_2}\right)^2 + 2 \omega \left(\alpha\dfrac{\partial f(z)}{\partial z_1} + \beta\dfrac{\partial f(z)}{\partial z_2}\right) f(z + c)  + b f(z + c)^2 = e^{g(z)}
\end{align*}
in $ \mathbb{C}^2 $ admits a transcendental entire solutions of finite order, then $ g(z) $ must be a polynomial in $ \mathbb{C}^2 $ of the form $ g(z)= L(z) + H(s) + B_1 $, where $ L(z)= a_1 z_1 + a_2 z_2 $ and $ H(s) $ is a polynomial in $ s:= c_2 z_1-c_1 z_2 $ in $ \mathbb{C}^2 $, $ a_1, a_2, B_1\in \mathbb{C} $. Further, $ f(z) $ must be one of the following forms:
\begin{enumerate}
	\item [(i)] 
	\begin{align*}
		f(z)= \pm \sqrt{\dfrac{(b\omega^2_2 + a)\omega^2_1 - 2(\omega_1\omega_2 b + a)\omega_1\omega_2 + (b\omega^2_1 + a)\omega^2_2} {ab\omega_1\omega_2(\omega_2 -\omega_1)^2}} e^{\frac{g(z-c)}{2}},
	\end{align*}
	where $ g(z)= \psi(z_2 - \frac{\beta}{\alpha}z_1) $ and $ \psi $ is a polynomial in $ \mathbb{C}^2 $;
	\item [(ii)] 
	\begin{align*}
		f(z) = \dfrac{\xi^2 -1}{\xi\sqrt{b}(\omega_2-\omega_1)} e^{\frac{L(z)+ H(s) + B}{2}},\; \mbox{where}\; B\in\mathbb{C}
	\end{align*}
 with 
	\begin{align*}
		\dfrac{\sqrt{a}(\xi^2 -1)}{2\sqrt{b}(\omega_2\xi^2 -\omega_1)}(\alpha a_{1} + \beta a_{2}) = e^{\frac{a_1 c_1 + a_2 c_2}{2}};
	\end{align*}
	\item [(iii)] 
	\begin{align*}
		f(z) = \dfrac{e^{L_1(z) +H_1(s) + E_1 -L_1(c)} - e^{L_2(z) +H_2(s) + E_2 -L_2(c)}}{\sqrt{b}(\omega_2 -\omega_1)},
	\end{align*}
	where $ L_{l}(z)= a_{l1} z_1  + a_{l2} z_2 $ for $ l=1, 2 $ with $ E_1, E_2\in\mathbb{C} $ and $ H_{l}(s) $ for $ l=1, 2 $ are polynomial in $ s:= c_2 z_1 -c_1z_2 $ in $ \mathbb{C}^2 $, satisfy
	\begin{align*}
		L_1(z) +H_1(s)\neq L_2(z) +H_2(s), \;\; g(z)= L_1(z) + L_2(z) + H_1(s) + H_2(s) + E_1 + E_2,
	\end{align*}
	and
	\begin{align*}
		\dfrac{\sqrt{a}}{\omega_2\sqrt{b}}[\alpha a_{11} +\beta a_{12}] e^{-L_1(c)}\equiv 1 \;\;\; \mbox{and} \;\;\; \dfrac{\sqrt{a}}{\omega_1\sqrt{b}}[\alpha a_{21} +\beta a_{22}] e^{-L_2(c)}\equiv 1.
	\end{align*}
\end{enumerate}	
\end{cor}
The following two examples validates the existence and form of the solutions of equations considered in Theorem \ref{th-2.1}.
\begin{example}
For $ c=( 7, -2, -4) $, by a routine computation, it can be easily shown that the transcendental entire solutions in $ \mathbb{C}^3 $ of the partial differential-difference equation	
\begin{align*}
	\left(2\dfrac{\partial f(z)}{\partial z_1} - \dfrac{\partial f(z)}{\partial z_3}\right)^2 -6 \left(2\dfrac{\partial f(z)}{\partial z_1} - \dfrac{\partial f(z)}{\partial z_3}\right) f(z + c)  &+ 2 f(z + c)^2 \\&= e^{4z_1 +\ln(6+6\sqrt{7}) z_2 +7z_3 +H(s) +\frac{\pi i}{3}},
\end{align*}
where $ H(s) $ is a polynomial in $ s= 2z_1 +z_2 +3z_3 $; must be of the form
\begin{align*}
	f(z_1,z_2,z_3) = \dfrac{1}{2\sqrt{14}} e^{\frac{1}{2}[4z_1 +\ln(6+6\sqrt{7}) z_2 +7z_3 +H(s) +\frac{\pi i}{3}]}.
\end{align*}
\end{example}

\begin{example}
For $ c=( 2, 3, 5) $, the transcendental entire solutions in $ \mathbb{C}^3 $ of the differential-difference equation	
\begin{align*}
	2\left(\dfrac{\partial f(z)}{\partial z_1} - 2\dfrac{\partial f(z)}{\partial z_3}\right)^2 &-8 \left(\dfrac{\partial f(z)}{\partial z_1} - 2\dfrac{\partial f(z)}{\partial z_3}\right) f(z + c)  + 3 f(z + c)^2 \\&= e^{15z_1 +\frac{1}{3}[\ln\left(\frac{9\sqrt{2}}{2\sqrt{2}\mp\sqrt{5}}\right) +\ln\left(\frac{18\sqrt{2}}{2\sqrt{2}\pm\sqrt{5}}\right)] z_2 -6z_3 +H(s) +\frac{16\pi i}{63}},
\end{align*}
where $ H(s) $ is a polynomial in $ s= 8z_1 + 3z_2 -5z_3 $; must be of the form	
\begin{align*}
	f(z_1,z_2,z_3) =& \dfrac{1}{\mp2\sqrt{5}}e^{5z_1 +\frac{1}{3}\ln\left(\frac{9\sqrt{2}} {2\sqrt{2}\mp\sqrt{5}}\right)z_2 -2z_3 +H_1(s) +\left(\frac{\pi i}{7} -\ln\left(\frac{9\sqrt{2}}{2\sqrt{2}\mp\sqrt{5}}\right)\right)} \\&- \dfrac{1}{\mp2\sqrt{5}} e^{10z_1 +\frac{1}{3}\ln\left(\frac{18\sqrt{2}} {2\sqrt{2}\pm\sqrt{5}}\right)z_2 -4z_3 +H_2(s) +\left(\frac{\pi i}{9} -\ln\left(\frac{18\sqrt{2}}{2\sqrt{2}\pm\sqrt{5}}\right)\right)}.
\end{align*}
\end{example}
The difference operator $ \Delta_cf $ of entire functions $ f $ in $ \mathbb{C}^n $ is defined by $ \Delta_cf(z):=f(z+c)-f(z) $. We obtain the following result finding the precise form of the solutions to a trinomial $ PDDE $ 
\begin{align}\label{eq-2.2}
	a\left(\alpha\dfrac{\partial f(z)}{\partial z_i} + \beta\dfrac{\partial f(z)}{\partial z_j}\right)^2 & + 2 \omega \left(\alpha\dfrac{\partial f(z)}{\partial z_i} + \beta\dfrac{\partial f(z)}{\partial z_j}\right) \Delta_cf(z) + b [\Delta_cf(z)]^2 = e^{g(z)},
\end{align}
involving $ \Delta_cf(z) $.
\begin{thm}\label{th-2.2}
Let $ c\in\mathbb{C}^n\setminus\{0\} $, $ a,b,\alpha\neq 0 $, $ \omega^2\neq ab $ and $ 1\leq i< j\leq n $, and $ \alpha d_{i} + \beta d_{j} \neq 0 $. If $ f(z) $ is a finite order transcendental entire solution of the $ PDDE $ \eqref{eq-2.2}, then $ f $ must assume one of the following forms:
\begin{enumerate}
	\item [(i)] 
	\begin{align*}
		f(z)= \phi\left(z_j - \frac{\beta}{\alpha}z_i\right),
	\end{align*} 
     where $ \phi $ is a finite order transcendental entire function satisfying
     \begin{align*}
     	\phi\left(z_j - \frac{\beta}{\alpha}z_i +c_j -\frac{\beta}{\alpha} c_i\right) - \phi\left(z_j - \frac{\beta}{\alpha}z_i\right)= M_2 e^{\frac{g(z-c)}{2}},
     \end{align*}
     
     \item [(ii)] 
     \begin{align*}
     	f(z) = \pm\frac{1}{\alpha\sqrt{a}}\int_{0}^{\frac{z_i}{\alpha}} e^{\frac{L(z) + H(s) + R}{2}} dz_i + \psi_1\left(z_j - \frac{\beta}{\alpha}z_i\right),
     \end{align*}
      $ g(z)= L(z) + H(s) + R $, where $ L(z)= a_1 z_1 +\cdots + a_n z_n $ and $ H(s) $ is a polynomial in $ s:= d_1 z_1 +\cdots + d_n z_n $ in $ \mathbb{C}^n $ with $ d_1 c_1 +\cdots + d_n c_n = 0 $ with $ H(z+c)= H(z) $, $ R\in\mathbb{C} $, and $ a_1 c_1 +\cdots + a_n c_n = 4k\pi i $ for $ k\in\mathbb{Z} $, and $ \psi_1 $ is a finite order periodic entire function with period $ (c_j -\frac{\beta}{\alpha} c_i) $; $ a_1,\ldots,a_n\in\mathbb{C} $.
      \item [(iii)] 
      \begin{align*}
      	f(z)= \dfrac{2(\omega_2\xi^2 -\omega_1)}{\xi\sqrt{a}(\omega_2-\omega_1)(k_i \alpha + k_j \beta)}e^{\frac{L(z) + R_2}{2}} + \phi_1\left(z_j -\frac{\beta}{\alpha} z_i\right),
      \end{align*}
      $ g(z)= L(z) + R_2 $, where $ L(z)=k_1 z_1 +\cdots + k_nz_n $, $ R_2 \in \mathbb{C} $, $ \phi_1 $ is a finite order periodic entire function with period $ (c_j -\frac{\beta}{\alpha} c_i) $ and satisfying
      \begin{align*}
      	\dfrac{\sqrt{a}(\xi^2 -1)}{2\sqrt{b}(\omega_2\xi^2 -\omega_1)}(\alpha k_i + \beta k_j) + 1 = e^{\frac{k_1 c_1 +\cdots + k_n c_n}{2}};
      \end{align*}
      \item [(iv)] 
      \begin{align*}
      	f(z) =\dfrac{1}{\sqrt{a} (\omega_2 -\omega_1)}\left[\dfrac{\omega_2 e^{L_1(z) + R_3}}{(\alpha a_{1i} + \beta a_{1j})} - \dfrac{\omega_1 e^{L_2(z) + R_4}}{(\alpha a_{2i} + \beta a_{2j})}\right] + \phi_2\left(z_j -\frac{\beta}{\alpha} z_i\right),
      \end{align*}
       $ g(z)=L_1(z) + L_2(z) + R_3 + R_4 $, $ L_1(z)\neq L_2(z) $ where $ L_{l}(z)= a_{l1} z_1 +\cdots + a_{ln} z_n $ and $ R_3, R_4\in\mathbb{C} $, $ \phi_2 $ is a finite order periodic function with period $ (c_j -\frac{\beta}{\alpha} c_i) $ and satisfying
      \begin{align*}
       \begin{cases}
       	\dfrac{\sqrt{a}}{\omega_2\sqrt{b}}[(\alpha a_{1i} +\beta a_{1j}) + \sqrt{a}\omega_2] e^{-L_1(c)}\equiv 1, \vspace{1.5mm}\\ \dfrac{\sqrt{a}}{\omega_1\sqrt{b}}[(\alpha a_{2i} +\beta a_{2j}) + \sqrt{b}\omega_1] e^{-L_2(c)}\equiv 1.
       \end{cases}
      \end{align*}
\end{enumerate}
\end{thm}
\begin{rem}
Theorem \ref{th-2.2} is a generalization of the binomial result of Xu \emph{et al.} \cite[Theorem 2.1]{Xu-ZHANG-ZHENG-RMJM-2022} in $ \mathbb{C}^n $.
\end{rem}
As a consequence of Theorem \ref{th-2.2}, we obtain the following corollary in $ \mathbb{C}^2 $ and it establishes solutions of trinomial $ PDDEs $ relating to a result concerning binomial $ PDDE $ \eqref{eq-1.3} in \cite{Xu-ZHANG-ZHENG-RMJM-2022}.
\begin{cor}
Let $ c\in\mathbb{C}^2\setminus\{0\} $, $ a,b,\alpha\neq 0 $, $ \omega^2\neq ab $ and $ \alpha d_{1} + \beta d_{2} \neq 0 $. Let $ f(z) $ be a finite order transcendental entire solution of the partial differential-difference equation
\begin{align*}
	a\left(\alpha\dfrac{\partial f(z)}{\partial z_1} + \beta\dfrac{\partial f(z)}{\partial z_2}\right)^2 & + 2 \omega \left(\alpha\dfrac{\partial f(z)}{\partial z_1} + \beta\dfrac{\partial f(z)}{\partial z_2}\right) \Delta_cf(z) \\&\nonumber + b [\Delta_cf(z)]^2 = e^{g(z)},
\end{align*}
then $ f(z) $ must satisfy one of the following cases:
\begin{enumerate}
	\item [(i)] 
	\begin{align*}
		f(z)= \phi(z_2 - \frac{\beta}{\alpha}z_1),
	\end{align*} 
	where $ \phi $ is a finite order transcendental entire function satisfying
	\begin{align*}
		\phi\left(z_2 - \frac{\beta}{\alpha}z_1 +c_2 -\frac{\beta}{\alpha} c_1\right) - \phi\left(z_2 - \frac{\beta}{\alpha}z_1\right)= \pm M_2 e^{\frac{g(z-c)}{2}}.
	\end{align*}
	\item [(ii)] 
	\begin{align*}
		f(z) = \pm\frac{1}{\alpha\sqrt{a}}\int_{0}^{\frac{z_1}{\alpha}} e^{\frac{L(z) + H(s) + R}{2}} dz_1 + \psi_1\left(z_2 - \frac{\beta}{\alpha}z_1\right),
	\end{align*}
	$ g(z)= L(z) + H(s) + R $, where $ L(z)= a_1 z_1 + a_2 z_2 $ and $ H(s) $ is a polynomial in $ s:= c_2 z_1-c_1 z_2 $ in $ \mathbb{C}^2 $, $ R\in\mathbb{C} $, and $ a_1 c_1 + a_2 c_2 = 4k\pi i $ for $ k\in\mathbb{Z} $, and $ \psi_1 $ is a finite order periodic entire function with period $ (c_2 -\frac{\beta}{\alpha} c_1) $.
	\item [(iii)] 
	\begin{align*}
		f(z)= \dfrac{2(\omega_2\xi^2 -\omega_1)} {\xi\sqrt{a}(\omega_2-\omega_1)(k_1 \alpha + k_2 \beta)}e^{\frac{k_1 z_1 + k_2 z_2 + R_2}{2}} + \phi_1(z_2 -\frac{\beta}{\alpha} z_1),
	\end{align*}
	$ g(z)= L(z) + R_2 $, where $ L(z)=k_1 z_1 + k_2 z_2 $, $ k_1, k_2, R_2 \in \mathbb{C} $; $ \phi_1 $ is a finite order periodic entire function with period $ (c_2 -\frac{\beta}{\alpha} c_1) $ and satisfying
	\begin{align*}
		\dfrac{\sqrt{a}(\xi^2 -1)}{2\sqrt{b}(\omega_2\xi^2 -\omega_1)}(\alpha k_1 + \beta k_2) + 1 = e^{\frac{k_1 c_1 + k_2 c_2}{2}};
	\end{align*}
	\item [(iv)] 
	\begin{align*}
		f(z) =\dfrac{1}{\sqrt{a} (\omega_2 -\omega_1)}\left(\dfrac{\omega_2 e^{L_1(z) + R_3}}{(\alpha a_{11} + \beta a_{12})} - \dfrac{\omega_1 e^{L_2(z) + R_4}}{(\alpha a_{21} + \beta a_{22})}\right) + \phi_2\left(z_2 -\frac{\beta}{\alpha} z_1\right),
	\end{align*}
	$ g(z)=L_1(z) + L_2(z) + R_3 + R_4 $, $ L_1(z)\neq L_2(z) $ where $ L_{l}(z)= a_{l1} z_1 + a_{l2} z_2 $ and $ R_3, R_4\in\mathbb{C} $, $ \phi_2 $ is a finite order periodic function with period $ (c_2 -\frac{\beta}{\alpha} c_1) $ and satisfying
	\begin{align*}
		\begin{cases}
			\dfrac{\sqrt{a}}{\omega_2\sqrt{b}}[(\alpha a_{11} +\beta a_{12}) + \sqrt{a}\omega_2] e^{-L_1(c)}\equiv 1, \vspace{1.5mm} \\ \dfrac{\sqrt{a}}{\omega_1\sqrt{b}}[(\alpha a_{21} +\beta a_{22}) + \sqrt{b}\omega_1] e^{-L_2(c)}\equiv 1.
		\end{cases}
	\end{align*}
\end{enumerate}
\end{cor}
The following examples are exhibited to validate the existence and precise form of the solutions of equations in Theorem \ref{th-2.2}.
\begin{example}
For $ c=( 2, 2, 3) $, the transcendental entire solutions in $ \mathbb{C}^3 $ of the differential-difference equation
\begin{align*}
\left(2\dfrac{\partial f(z)}{\partial z_1} + \dfrac{\partial f(z)}{\partial z_3}\right)^2  - 8 \left(2\dfrac{\partial f(z)}{\partial z_1} + \dfrac{\partial f(z)}{\partial z_3}\right) \Delta_cf(z) + 3 [\Delta_cf(z)]^2 = e^{3z_1 +\ln\left(\frac{6+3\sqrt{13}}{4 +\sqrt{13}}\right) z_2 -2z_3 + \frac{\pi i}{7}},
\end{align*}
must be of the form
\begin{align*}
	f(z_1,z_2,z_3)=\dfrac{\sqrt{3}(4 +3\sqrt{13})}{4\sqrt{26}}e^{\frac{1}{2}[3z_1 +\ln\left(\frac{6+3\sqrt{13}}{4 +\sqrt{13}}\right) z_2 -2z_3 + \frac{\pi i}{7}]} + e^{\pi i\left(\frac{z_1}{2} + z_3\right)}.
\end{align*}
\end{example}

\begin{example}
For $ c=( 3, 1, -4) $, the transcendental entire solutions in $ \mathbb{C}^3 $ of the differential-difference equation	
\begin{align*}
	3\left(3\dfrac{\partial f(z)}{\partial z_1} + 2\dfrac{\partial f(z)}{\partial z_3}\right)^2 & -10 \left(3\dfrac{\partial f(z)}{\partial z_1} + 2\dfrac{\partial f(z)}{\partial z_3}\right) \Delta_cf(z) + [\Delta_cf(z)]^2 \\&= e^{12z_1 +\left(\ln\left(\frac{3(23\mp\sqrt{22})}{5\mp\sqrt{22}}\right)+ \ln\left(\frac{\sqrt{3}(36\sqrt{3} +5 \pm\sqrt{22})}{5\pm\sqrt{22}} \right) \right)z_2 + 9z_3 + \frac{(2\pi i +\sqrt{5} + \sqrt{3})}{\sqrt{7}}},
\end{align*}	
must be of the form
\begin{align*}
	f(z_1,z_2,z_3) &=\dfrac{(5\mp \sqrt{22})}{\mp36\sqrt{66}}e^{4z_1 +\ln\left(\frac{3(23\mp\sqrt{22})}{5\mp\sqrt{22}}\right)z_2 +3z_3 + \frac{(\pi i +\sqrt{3})}{\sqrt{7}}} \\&- \dfrac{(5\pm \sqrt{22})}{\mp72\sqrt{66}}e^{8z_1 +\ln\left(\frac{\sqrt{3}(36\sqrt{3} +5\pm\sqrt{22})}{5\pm\sqrt{22}}\right)z_2 +6z_3 + \frac{(\pi i+\sqrt{5})}{\sqrt{7}}} + e^{\pi i\left(\frac{2z_1}{3} +z_3\right)}.
\end{align*}
\end{example}
 \section{Key lemmas and Proof of the main results}
 First, we present here some necessary lemmas which will play a key roles in proving the main results of this paper.
 	
\begin{lem}\cite{Ronkin_AMS-1974,Stoll-AMS-1974}\label{lem-3.1}
For any entire function $ F $ on $ \mathbb{C}^n $, $ F(0)\neq 0 $ and put $\rho(n_F)=\rho < \infty $, where $ \rho(n_F) $ denotes be the order of the counting function of zeros of $ F $. Then there exist a canonical function $ f_F $ and a function $ g_F \in\mathbb{C}^n $ such that $ F(z) =f_F (z)e^{g_F (z)} $. For the special case $ n = 1 $, $ f_F $ is the canonical product of Weierstrass.
 \end{lem}
 \begin{lem}\cite{Polya-JLMS-1926}\label{lem-3.2}
 If $  g $ and $ h $ are entire functions on the complex plane $ \mathbb{C} $ and $ g(h) $ is an entire function of finite order, then there are only two possible cases: either
 \begin{enumerate}
 	\item [(i)] the internal function $ h $ is a polynomial and the external function $ g $ is of finite order; or
   \item[(ii)] the internal function $ h $ is not a polynomial but a function of finite order, and the external function $ g $ is of zero order.
 \end{enumerate}
 \end{lem}
\begin{lem}\cite{Hu-Li-Yang-2003}\label{lem-3.3}
Suppose that $ a_0(z), a_1(z),\ldots, a_m(z) \; (m\geq 1) $ are meromorphic functions on $ \mathbb{C}^n $ and $ g_0(z), g_1(z),\ldots, g_m(z) $ are entire functions on $ \mathbb{C}^n $ such that $ g_i(z)- g_j(z) $ are not constants for $ 0\leq i< j\leq m $. If
\begin{align*}
	\sum_{i=0}^{m} a_i(z) e^{g_i(z)} \equiv 0
\end{align*}
and $ || T(r, a_i) = o(T(r)), \;\;\;\; i=0, 1, \ldots, m $ holds, where $ T(r):=\min_{0\leq i< j\leq m} T(r, e^{g_i -g_j}) $, then $ a_i(z) \equiv 0 $ for $ i=0, 1, \ldots, m $.
\end{lem}

\begin{lem}\cite{Hu-Li-Yang-2003}\label{lem-3.4}
 Let $ f_j(\not \equiv 0) $, $ j=1,2,3 $, be meromorphic functions on $ \mathbb{C}^n $ such that $ f_1 $ is non-constant and $ f_1+f_2+f_3=1 $ such that 
\begin{align*}
	\sum_{j=1}^{3}\left\{N_2(r,\frac{1}{f_j})+2\overline{N}(r,f_j)\right\}<\lambda T(r,f_1) + O(\log^{+} T(r,f_1)),
\end{align*}
for all $ r $ outside possibly a set with finite logarithmic measure, where $ \lambda<1 $ is a positive number. Then either $ f_2=1\;\mbox{or}\; f_3=1 $.
\end{lem}
\begin{rem}
	Here, $ N_2(r,1/f) $ is the counting function of the zeros of $ f $ in $ |z|\leq r $, where the	simple zero is counted once, and the multiple zero is counted twice.
\end{rem}

Now we discuss the proof of the main results of the paper.

\begin{proof}[\bf Proof of Theorem \ref{th-2.1}]
Assume that $ f $ is a transcendental entire solution of finite order of the equation \eqref{eq-2.1}. we see that \eqref{eq-2.1} can be written as
\begin{align}\label{eq-3.1}
	aF^2 + 2\omega F G + bG^2 =1,
\end{align}
where, $ F $ and $ G $ are defined by 
\begin{align}\label{eq-3.2}
	F := \dfrac{\alpha\dfrac{\partial f(z)}{\partial z_i} + \beta\dfrac{\partial f(z)}{\partial z_j}}{e^{\frac{g(z)}{2}}} \;\;\;\; \mbox{and} \;\;\;\; G:= \dfrac{f(z + c)} {e^{\frac{g(z)}{2}}}.
\end{align}
It is easy to see that \eqref{eq-3.1} can be expressed as
\begin{align*}
	(\sqrt{a} F -\omega_1\sqrt{b} G)(\sqrt{a} F -\omega_2\sqrt{b} G) = 1,
\end{align*}
where $ \omega_1 =-\frac{\omega}{\sqrt{ab}} \pm \frac{\sqrt{\omega^2 -ab}}{\sqrt{ab}} $ and $ \omega_2 =-\frac{\omega}{\sqrt{ab}} \mp \frac{\sqrt{\omega^2 -ab}}{\sqrt{ab}} $. Since $ f $ is a finite order transcendental entire function and $ g $ is a polynomial, by Lemmas \ref{lem-3.1} and \ref{lem-3.2}, there exists a polynomial $ p $ in $ \mathbb{C}^n $ such that
\begin{align}\label{eq-3.3}
	\sqrt{a} F -\omega_1\sqrt{b} G = e^{p} \;\;\;\;\mbox{and} \;\;\;\; \sqrt{a} F -\omega_2\sqrt{b} G =e^{-p}.
\end{align}
An elementary computation using \eqref{eq-3.2} and \eqref{eq-3.3} shows that
\begin{align}\label{eq-3.4}
	\alpha\dfrac{\partial f(z)}{\partial z_i} + \beta\dfrac{\partial f(z)}{\partial z_j} = \dfrac{\omega_2 e^{p(z)} - \omega_1 e^{-p(z)}}{\sqrt{a}(\omega_2 -\omega_1)} e^{\frac{g(z)}{2}}, 
\end{align}
\begin{align}\label{eq-3.5}
	f(z+c) = \dfrac{e^{p(z)} - e^{-p(z)}}{\sqrt{b}(\omega_2 -\omega_1)} e^{\frac{g(z)}{2}}.
\end{align}
For brevity, we assume that 
\begin{align}\label{eq-3.6}
	h_1(z)= \dfrac{g(z)}{2} + p(z) \;\;\;\;\mbox{and}\;\;\;\; h_2(z)= \dfrac{g(z)}{2} - p(z).
\end{align} 
Therefore, the equations \eqref{eq-3.4} and \eqref{eq-3.5} can be written as
\begin{align}\label{eq-3.7}
	\alpha\dfrac{\partial f(z)}{\partial z_i} + \beta\dfrac{\partial f(z)}{\partial z_j} = \dfrac{\omega_2 e^{h_1(z)} - \omega_1 e^{h_2(z)}}{\sqrt{a}(\omega_2 -\omega_1)},
\end{align}
\begin{align}\label{eq-3.8}
	f(z+c) = \dfrac{e^{h_1(z)} - e^{h_2(z)}}{\sqrt{b}(\omega_2 -\omega_1)}.
\end{align}
In view of \eqref{eq-3.7} and \eqref{eq-3.8}, a simple computation yields that
\begin{align}\label{eq-3.9}
	H_{11}(z) e^{h_1(z) -h_1(z+c)} - H_{12}(z) e^{h_2(z) -h_1(z+c)} + K_{1} e^{h_2(z+c) -h_1(z+c)} \equiv 1,
\end{align}
where 
\begin{align*}
	\begin{cases}
		H_{11}(z) = \dfrac{\sqrt{a}\left(\alpha\dfrac{\partial h_1(z)}{\partial z_i} + \beta\dfrac{\partial h_1(z)}{\partial z_j}\right)} {\omega_2\sqrt{b}},\vspace{1.2mm}\\
		H_{12}(z) = \dfrac{\sqrt{a}\left(\alpha\dfrac{\partial h_2(z)}{\partial z_i} + \beta\dfrac{\partial h_2(z)}{\partial z_j}\right)} {\omega_2\sqrt{b}},\; K_{1} = \dfrac{\omega_1}{\omega_2}.
	\end{cases}
\end{align*}
\noindent {\bf Case A:} If $ e^{h_2(z+c) -h_1(z+c)} $ is a constant, then $ h_2(z+c)-h_1(z+c) =K $, where $ K\in\mathbb{C}^n $ is a constant. From \eqref{eq-3.6}, it is easy to see that $ p(z)= -K $ is a constant. Let $ \xi =e^{p(z)} $, then the equations \eqref{eq-3.7} and \eqref{eq-3.8} become
\begin{align}\label{eq-3.10}
	\alpha\dfrac{\partial f(z)}{\partial z_i} + \beta\dfrac{\partial f(z)}{\partial z_j} = M_{1} e^{\frac{g(z)}{2}} \;\; \mbox{and}\;\; f(z+c)= M_{2} e^{\frac{g(z)}{2}}
\end{align}
where 
\begin{align*}
	M_{1} = \dfrac{\omega_2 \xi - \omega_1 \xi^{-1}}{\sqrt{a}(\omega_2 -\omega_1)}, \;\;\;\; M_{2}= \dfrac{\xi - \xi^{-1}}{\sqrt{b}(\omega_2 -\omega_1)}
\end{align*}
and
\begin{align}\label{eq-3.11}
	M^2_{1} + M^2_{2} = \dfrac{(b\omega^2_2 + a)\xi^2 - 2(\omega_1\omega_2 b + a) + (b\omega^2_1 + a)\frac{1}{\xi^2}}{ab(\omega_2 -\omega_1)^2}.
\end{align}
It is easy to see that $ M_{2}\neq 0 $. Now, we discuss the following two sub-cases. \vspace{1.2mm}

\noindent{\bf Sub-case A1:} Suppose that $ M_{1}= 0 $, clearly $ \xi^2 ={\omega_1}/{\omega_2} $. In view of \eqref{eq-3.11}, it is easy to see that
\begin{align*}
	M_{2} = \pm \sqrt{\dfrac{(b\omega^2_2 + a)\omega^2_1 - 2(\omega_1\omega_2 b + a)\omega_1\omega_2 + (b\omega^2_1 + a)\omega^2_2} {ab\omega_1\omega_2(\omega_2 -\omega_1)^2}}\;(=N_1; \;\mbox{says}).
\end{align*}
Thus, the equation \eqref{eq-3.10} becomes
\begin{align}\label{eq-3.12}
	&\alpha\dfrac{\partial f(z)}{\partial z_i} + \beta\dfrac{\partial f(z)}{\partial z_j} = 0 \;\; \mbox{and}\;\;\\&\label{eq-3.13} f(z+c)= N_1 e^{\frac{g(z)}{2}}.
\end{align}
Solving \eqref{eq-3.12}, we obtain
\begin{align}\label{eq-3.14}
	f(z)= \phi\left(z_j - \frac{\beta}{\alpha}z_i\right),
\end{align}
where $ \phi\left(z_j - \frac{\beta}{\alpha}z_i\right) $ is a finite order transcendental entire function. Furthermore, the equation \eqref{eq-3.13} can be written as
\begin{align}\label{eq-3.15}
	f(z)= \pm \sqrt{\dfrac{(b\omega^2_2 + a)\omega^2_1 - 2(\omega_1\omega_2 b + a)\omega_1\omega_2 + (b\omega^2_1 + a)\omega^2_2} {ab\omega_1\omega_2(\omega_2 -\omega_1)^2}} e^{\frac{g(z-c)}{2}}.
\end{align}
In view of \eqref{eq-3.14} and \eqref{eq-3.15}, we obtain
\begin{align*}
	g(z)= \psi\left(z_j - \frac{\beta}{\alpha}z_i\right)= 2\ln\left(\pm\dfrac{\phi(z_j - \frac{\beta}{\alpha}z_i +c_j -\frac{\beta}{\alpha} c_i)(\omega_2 -\omega_1)\sqrt{ab\omega_1\omega_2}} {\sqrt{(b\omega^2_2 + a)\omega^2_1 - 2(\omega_1\omega_2 b + a)\omega_1\omega_2 + (b\omega^2_1 + a)\omega^2_2}}\right).
\end{align*}
\noindent{\bf Sub-case A2:} Suppose that $ M_{1}\neq 0 $. Then, it follows from \eqref{eq-3.10} we obtain
\begin{align}\label{eq-3.16}
	\dfrac{M_{2}}{2M_{1}}\left(\alpha\dfrac{\partial g(z)}{\partial z_i} + \beta\dfrac{\partial g(z)}{\partial z_j}\right) = e^{\frac{g(z+c) -g(z)}{2}}.
\end{align}
Since $ g(z) $ is a polynomial, \eqref{eq-3.16} implies that $ g(z+c) -g(z) =\xi_1 $, where $ \xi_1 $ is a constant in $ \mathbb{C} $. Therefore, it follows that $ g(z)= L(z) + H(s) + B_1 $, where $ L(z)= a_1 z_1 +\cdots + a_n z_n $ and $ H(s) $ is a polynomial in $ s:= d_1 z_1 +\cdots + d_n z_n $ in $ \mathbb{C}^n $ with $ d_1 c_1 +\cdots + d_n c_n = 0 $ with $ H(z+c)= H(z) $. Thus from \eqref{eq-3.16}, we obtain
\begin{align*}
	\alpha\dfrac{\partial L(z)}{\partial z_i} + \beta\dfrac{\partial L(z)}{\partial z_j} + \alpha\dfrac{\partial H(s)}{\partial z_i} + \beta\dfrac{\partial H(s)}{\partial z_j} \equiv M_{3},
\end{align*}
or,
\begin{align*}
	\alpha\dfrac{\partial H(s)}{\partial z_i} + \beta\dfrac{\partial H(s)}{\partial z_j}\equiv (\alpha d_{i} + \beta d_{j}) H^{\prime}\equiv M_{4},
\end{align*}
where $ M_{3}=\dfrac{2M_{1}}{M_{2}} e^{\frac{a_1 c_1 +\cdots + a_n c_n}{2}} $ and $ M_{4} =M_{3} - (\alpha a_{i} + \beta a_{j}) $.\vspace{1.2mm}

If $ \alpha d_{i} + \beta d_{j} \neq 0 $, then $ H^{\prime} $ is a constant in $ \mathbb{C}^n $. If follows that $ H(s)= A_1 s + A_2 = A_1(d_1 z_1 +\cdots + d_n z_n) + A_2 $, where $ A_1 =\frac{M_4}{\alpha d_{i} + \beta d_{j}}$ and $ A_2 \in\mathbb{C} $. Therefore, $ L(z) + H(s) $ is also a linear function. For convenience, we still denote $ g(z)= L(z) + B_{1} $, which implies that $ H(s)=0 $. Thus, we obtain $ \xi_1 = a_1 c_1 +\cdots + a_n c_n $. \vspace{1.2mm}

If $ \alpha d_{i} + \beta d_{j} = 0 $, then $ (\alpha d_{i} + \beta d_{j}) H^{\prime}\equiv 0 $. It follows that, $ M_{4} =M_{3} - (\alpha a_{i} + \beta a_{j}) \equiv 0 $.\vspace{1.2mm}

A simple computation shows that
\begin{align*}
	\dfrac{\sqrt{a}(\xi^2 -1)}{2\sqrt{b}(\omega_2\xi^2 -\omega_1)}(\alpha a_{i} + \beta a_{j}) = e^{\frac{a_1 c_1 +\cdots + a_n c_n}{2}}.
\end{align*}
Hence, from the second equation of \eqref{eq-3.10}, we obtain
\begin{align*}
	f(z)= M_{2} e^{\frac{g(z-c)}{2}} = \dfrac{\xi^2 -1}{\xi\sqrt{b}(\omega_2-\omega_1)} e^{\frac{L(z)+ H(s) + B}{2}}, \;\;\mbox{where}\;\; B= B_1 -L(c).
\end{align*}

\noindent{\bf Case B:} If $ e^{h_2(z+c) -h_1(z+c)} $ is not a constant, then obviously, $ H_{11}(z)\equiv 0 $ and $ H_{12}(z)\equiv 0 $ cannot hold at the same time. Otherwise, from \eqref{eq-3.9} we see that $ K_{1} e^{h_2(z+c) -h_1(z+c)}\equiv 1 $, a contradiction. If $ H_{11}(z)\equiv 0 $ and $ H_{12}(z)\not\equiv 0 $, then in view of \eqref{eq-3.9}, we obtain
\begin{align}\label{eq-3.17}
	- H_{12}(z) e^{h_2(z) -h_1(z+c)} + K_{1} e^{h_2(z+c) -h_1(z+c)} \equiv 1.
\end{align}
Since $ e^{h_2(z+c) -h_1(z+c)} $ is not a constant, it follows that $ e^{h_2(z) -h_1(z+c)} $ is not a constant. Furthermore, $ e^{h_2(z+c) -h_2(z)} $ is not a constant. Otherwise, $ h_2(z+c) -h_2(z)=\xi_2 $, where $ \xi_2\in\mathbb{C} $. Then, from \eqref{eq-3.17} we see that $ (-H_{12}(z) e^{-\xi_2} + K_1) e^{h_2(z+c) -h_1(z+c)} \equiv 1 $, which is a contradiction as $ e^{h_2(z+c) -h_1(z+c)} $ is non-constant. Therefore, the equation \eqref{eq-3.17} can be expressed as 
\begin{align}\label{eq-3.18}
	- H_{12}(z) e^{h_2(z)} + K_{1} e^{h_2(z+c)} - e^{h_1(z+c)} \equiv 0.
\end{align}
In view of Lemma \ref{lem-3.3}, from \eqref{eq-3.18}, we get a contradiction. Similarly, if $ H_{11}(z)\not\equiv 0 $ and $ H_{12}(z)\equiv 0 $, we can get a contradiction. Thus, we conclude that $ H_{11}(z)\not\equiv 0 $ and $ H_{12}(z)\not\equiv 0 $.

As $ h_1(z), h_2(z) $ are polynomials and $ K_{1} e^{h_2(z+c) -h_1(z+c)} $ is non-constant, then by Lemma \ref{lem-3.4} for \eqref{eq-3.9}, we obtain
\begin{align}\label{eq-3.19}
	H_{11}(z) e^{h_1(z) -h_1(z+c)}\equiv 1 \;\;\mbox{or}\;\; - H_{12}(z) e^{h_2(z) -h_1(z+c)}\equiv 1.
\end{align}
\noindent{\bf Sub-case B1:} Assume that $ H_{11}(z) e^{h_1(z) -h_1(z+c)}\equiv 1 $. Then, from \eqref{eq-3.9} it is easy to see that $ \frac{H_{12}(z)}{K_{1}} e^{h_2(z) -h_2(z+c)}\equiv 1 $. Since $ h_1(z), h_2(z) $ are polynomials, it follows that $ h_1(z) -h_1(z+c)=\xi_3 $ and $ h_2(z) -h_2(z+c)=\xi_4 $, where $ \xi_3, \xi_4\in\mathbb{C} $. Thus, it follows that $ h_1(z)=L_1(z) +H_1(s) + E_1 $ and $ h_2(z)=L_2(z) +H_2(s) + E_2 $, where $ L_{l}(z)= a_{l1} z_1 +\cdots + a_{ln} z_n $ and $ H_{l}(s) $ for $ l=1, 2 $ are polynomial in $ s:= d_1 z_1 +\cdots + d_n z_n $ in $ \mathbb{C}^n $ with $ d_1 c_1 +\cdots + d_n c_n = 0 $ with $ H_{l}(z+c)= H_{l}(z) $ for $ l=1, 2 $, and $ E_1, E_2\in\mathbb{C} $. Obviously $ L_1(z) +H_1(s)\neq L_2(z) +H_2(s) $. Otherwise, $ h_2(z+c) -h_1(z+c) $ is a constant, which shows that $ e^{h_2(z+c) -h_1(z+c)} $ is a constant, a contradiction. Substituting $ h_1(z) $ and $ h_2(z) $ into $ H_{11}(z) e^{h_1(z) -h_1(z+c)}\equiv 1 $ and $ \frac{H_{12}(z)}{K_{1}} e^{h_2(z) -h_2(z+c)}\equiv 1 $, we obtain
\[
\begin{cases}
\dfrac{\sqrt{a}}{\omega_2\sqrt{b}}[(\alpha a_{1i} +\beta a_{1j}) + (\alpha d_i +\beta d_j)H^{\prime}_{1}] e^{-L_1(c)}\equiv 1,\vspace{1.3mm}\\
\dfrac{\sqrt{a}}{\omega_1\sqrt{b}}[(\alpha a_{2i} +\beta a_{2j}) + (\alpha d_i +\beta d_j)H^{\prime}_{2}] e^{-L_2(c)}\equiv 1.
\end{cases}
\]
By the similar argument used in Case A, we easily obtain $ (\alpha d_i +\beta d_j)H^{\prime}_{1}\equiv 0 $ and $ (\alpha d_i +\beta d_j)H^{\prime}_{2} \equiv 0 $, which implies that
\begin{align*}
	\begin{cases}
	  \dfrac{\sqrt{a}}{\omega_2\sqrt{b}}[\alpha a_{1i} +\beta a_{1j}] e^{-L_1(c)}\equiv 1 ,\vspace{1.3mm}\\ \dfrac{\sqrt{a}}{\omega_1\sqrt{b}}[\alpha a_{2i} +\beta a_{2j}] e^{-L_2(c)}\equiv 1.
	\end{cases}
\end{align*}
Therefore, from \eqref{eq-3.8}, we see that
\begin{align*}
	f(z) = \dfrac{e^{L_1(z) +H_1(s) + E_1 -L_1(c)} - e^{L_2(z) +H_2(s) + E_2 -L_2(c)}}{\sqrt{b}(\omega_2 -\omega_1)}
\end{align*}
From \eqref{eq-3.6}, it is easy to see that
\begin{align*}
	g(z)= h_1(z) + h_2(z) = L(z)+ H(s) + E,
\end{align*}
where $ L(z)=L_1(z) + L_2(z) $, $ H(s)= H_1(s) + H_2(s) $ and $ E =E_1 + E_2 $.
\vspace{1.2mm}

\noindent{\bf Sub-case B2:} Assume that $ - H_{12}(z) e^{h_2(z) -h_1(z+c)}  \equiv 1 $. Then, from \eqref{eq-3.9} it is easy to see that $ -\frac{H_{11}(z)}{K_{1}} e^{h_1(z) -h_2(z+c)}\equiv 1 $. Since $ h_1(z)$ and $ h_2(z) $ are polynomials, it follows that $ h_2(z) -h_1(z+c)=\xi_5 $ and $ h_1(z) -h_2(z+c)=\xi_6 $, where $ \xi_5, \xi_6\in\mathbb{C} $. A simple computation shows that $ h_1(z+2c) -h_1(z)= -\xi_5 -\xi_6 $ and $ h_2(z+2c) -h_2(z)= -\xi_5 -\xi_6 $. Thus, we deduce that $ h_1(z)= L(z) + H(s) + E_3 $ and $ h_2(z)= L(z) + H(s) + E_4 $, where $ L(z)= a_1 z_1 +\cdots + a_n z_n $ and $ H(s) $ is a polynomial in $ s:= d_1 z_1 +\cdots + d_n z_n $ in $ \mathbb{C}^n $ with $ d_1 c_1 +\cdots + d_n c_n = 0 $ with $ H(z+c)= H(z) $, and $ E_3, E_4\in \mathbb{C} $. Now, we see that $ h_2(z+c) -h_1(z+c)= E_4 - E_3 $, which shows that $ e^{h_2(z+c) -h_1(z+c)} $ is a constant, a contradiction. This completes the proof.
\end{proof}	

\begin{proof}[\bf Proof of Theorem \ref{th-2.2}]
Suppose that $ f(z) $ is a finite order transcendental entire solution of \eqref{eq-2.2}. The equation \eqref{eq-2.2} can be written as
\begin{align}\label{eq-3.20}
	(\sqrt{a} F -\omega_1\sqrt{b} G)(\sqrt{a} F -\omega_2\sqrt{b} G) = 1,
\end{align}
where $ F $ and $ G $ are defined by 
\begin{align}\label{eq-3.21}
	F := \dfrac{\alpha\dfrac{\partial f(z)}{\partial z_i} + \beta\dfrac{\partial f(z)}{\partial z_j}}{e^{\frac{g(z)}{2}}} \;\;\;\; \mbox{and} \;\;\;\; G:= \dfrac{f(z + c) -f(z)} {e^{\frac{g(z)}{2}}}.
\end{align}
By the similar argument being used in the proof of the Theorem \ref{th-2.1}, there exists a polynomial $ p $ in $ \mathbb{C}^n $ such that
\begin{align}\label{eq-3.22}
	\sqrt{a} F -\omega_1\sqrt{b} G = e^{p} \;\;\;\;\mbox{and} \;\;\;\; \sqrt{a} F -\omega_2\sqrt{b} G =e^{-p}.
\end{align}
A simple computation using \eqref{eq-3.21} and \eqref{eq-3.22} given us
\begin{align}\label{eq-3.23}
	\alpha\dfrac{\partial f(z)}{\partial z_i} + \beta\dfrac{\partial f(z)}{\partial z_j} = \dfrac{\omega_2 e^{h_1(z)} - \omega_1 e^{h_2(z)}}{\sqrt{a}(\omega_2 -\omega_1)}
\end{align}
\begin{align}\label{eq-3.24}
   f(z + c)-f(z)= \dfrac{e^{h_1(z)} - e^{h_2(z)}}{\sqrt{b}(\omega_2 -\omega_1)},
\end{align}	
where 
\begin{align}\label{eq-3.25}
	h_1(z)= \dfrac{g(z)}{2} + p(z) \;\;\;\;\mbox{and}\;\;\;\; h_2(z)= \dfrac{g(z)}{2} - p(z).
\end{align}
Thus, it follows from \eqref{eq-3.23} and \eqref{eq-3.24} that
\begin{align}\label{eq-3.26}
	H_{21}(z) e^{h_1(z) -h_1(z+c)} - H_{22}(z) e^{h_2(z) -h_1(z+c)} + K_{2} e^{h_2(z+c) -h_1(z+c)} \equiv 1,
\end{align}
where 
\begin{align*}
	\begin{cases}
		H_{21}(z) = \dfrac{\sqrt{a}\left(\alpha\dfrac{\partial h_1(z)}{\partial z_i} + \beta\dfrac{\partial h_1(z)}{\partial z_j}\right) + \sqrt{b}\omega_2} {\omega_2\sqrt{b}},\vspace{2mm}\\
		H_{22}(z) = \dfrac{\sqrt{a}\left(\alpha\dfrac{\partial h_2(z)}{\partial z_i} + \beta\dfrac{\partial h_2(z)}{\partial z_j}\right)+ \sqrt{b}\omega_1} {\omega_2\sqrt{b}}\;\;\mbox{and}\;\; K_2=\dfrac{\omega_1}{\omega_2}.
	\end{cases}
\end{align*}
\noindent{\bf Case A:} If $ e^{h_2(z+c) -h_1(z+c)} $ is a constant, then $ h_2(z+c) -h_1(z+c) =K $, where $ K\in\mathbb{C}^n $ is a constant. From \eqref{eq-3.25}, it is easy to see that $ p(z)= -K $ is a constant. Let $ \xi =e^{p(z)} $, then the equations \eqref{eq-3.23} and \eqref{eq-3.24} becomes
\begin{align}\label{eq-3.27}
	\alpha\dfrac{\partial f(z)}{\partial z_i} + \beta\dfrac{\partial f(z)}{\partial z_j} = M_{1} e^{\frac{g(z)}{2}} \;\; \mbox{and}\;\; f(z + c)-f(z)= M_{2} e^{\frac{g(z)}{2}}
\end{align}
where $ M_{1} $ and $ M_{2} $ are same as in Case A in the proof of Theorem \ref{th-2.1}. \vspace{1.2mm}

\noindent{\bf Sub-case A1:} Assume that $ M_{1}= 0 $, then we obtain $ \xi^2 ={\omega_1}/{\omega_2} $. By the similar argument being used in the proof of Theorem \ref{th-2.1}, we see that
\begin{align*}
	f(z)= \phi(z_j - \frac{\beta}{\alpha}z_i),
\end{align*}
where $ \phi(z_j - \frac{\beta}{\alpha}z_i) $ is a finite order transcendental entire function satisfying
\begin{align*}
	\phi\left(z_j - \frac{\beta}{\alpha}z_i +c_j -\frac{\beta}{\alpha} c_i\right) - \phi\left(z_j - \frac{\beta}{\alpha}z_i\right)= N_1 e^{\frac{g(z-c)}{2}},
\end{align*}
where $ N_1 $ is defined in Case A in the proof of Theorem \ref{th-2.1}. \vspace{1.2mm}

\noindent{\bf Sub-case A2:} If $ M_2 =0 $, then we see that $ \xi^2 =1 $. Using \eqref{eq-3.11} a simple computation shows that $ M_1 = \pm\frac{1}{\sqrt{a}} $. Therefore, from \eqref{eq-3.27} it follows that
\begin{align}\label{eq-3.28}
	\alpha\dfrac{\partial f(z)}{\partial z_i} + \beta\dfrac{\partial f(z)}{\partial z_j} = \pm\frac{1}{\sqrt{a}} e^{\frac{g(z)}{2}} \;\; \mbox{and}\;\; f(z + c) = f(z).
\end{align}
We see that
\begin{align*}
	&\alpha\dfrac{\partial f(z + c)}{\partial z_i} + \beta\dfrac{\partial f(z + c)}{\partial z_j} = \alpha\dfrac{\partial f(z)}{\partial z_i} + \beta\dfrac{\partial f(z)}{\partial z_j},
\end{align*}
which implies that $ e^{\frac{g(z+c)-g(z)}{2}} =1 $. Thus, we have $ g(z)= L(z) + H(s) + R $, where $ L(z)= a_1 z_1 +\cdots + a_n z_n $ and $ H(s) $ is a polynomial in $ s:= d_1 z_1 +\cdots + d_n z_n $ in $ \mathbb{C}^n $ with $ d_1 c_1 +\cdots + d_n c_n = 0 $ with $ H(z+c)= H(z) $, $ R\in\mathbb{C} $, and $ a_1 c_1 +\cdots + a_n c_n = 4k\pi i $ for $ k\in\mathbb{Z} $. The characteristic
equations for the first equation of \eqref{eq-3.28} are
\begin{align*}
	\dfrac{dz_i}{dt}= \alpha, \;\; \dfrac{dz_j}{dt}= \beta, \;\; \dfrac{df}{dt} =\pm\frac{1}{\sqrt{a}} e^{\frac{g(z)}{2}}.
\end{align*}
Using the initial conditions: $ z_i =0 $, $ z_j =s $ and $ f=f(0,s):= \psi_1(s) $, with a parameter $ s $. Therefore, we obtain the following parametric representation for the solutions of the characteristic equations: $ z_i=\alpha t $, $ z_j=\beta t + s $,
\begin{align*}
	f(s,t) = \pm\frac{1}{\sqrt{a}}\int_{0}^{t} e^{\frac{g(z)}{2}} dt + \psi_1(s)
\end{align*}
or,
\begin{align*}
	f(z) = \pm\frac{1}{\alpha\sqrt{a}}\int_{0}^{\frac{z_i}{\alpha}} e^{\frac{a_1 z_1 +\cdots + a_n z_n + H(d_1 z_1 +\cdots + d_n z_n) + R}{2}} dz_i + \psi_1\left(z_j - \frac{\beta}{\alpha}z_i\right),
\end{align*}
where, $ \psi_1 $ is a finite order entire function. Substituting $ f(z) $ into the second equation of \eqref{eq-3.28}, we obtain
\begin{align*}
\psi_1\left(z_j - \frac{\beta}{\alpha}z_i +c_j -\frac{\beta}{\alpha} c_i\right) = \psi_1\left(z_j - \frac{\beta}{\alpha}z_i\right),
\end{align*}
which implies that $ \psi_1 $ is a periodic function with period $ (c_j -\frac{\beta}{\alpha} c_i) $. \vspace{1.2mm}

\noindent{\bf Sub-case A3:} Suppose that $ M_1\neq 0 $ and $ M_2\neq 0 $. Then, from \eqref{eq-3.27} a simple computation shows that
\begin{align}\label{eq-3.29}
	\dfrac{M_2}{2M_1}\left(\alpha\dfrac{\partial g(z)}{\partial z_i} + \beta\dfrac{\partial g(z)}{\partial z_j}\right) + 1 = e^{\frac{g(z+c)-g(z)}{2}}.
\end{align}
As $ g(z) $ is a polynomial, from \eqref{eq-3.29} it follows that $ g(z+c)-g(z) =\eta $, where $ \eta $ is a constant in $ \mathbb{C} $. It yields that $ g(z)= L_1(z) + H(s)+ R_1 $, where $ L_1(z)= a_{11} z_1 +\cdots + a_{1n} z_n $ and $ H(s) $ is a polynomial in $ s:= d_1 z_1 +\cdots + d_n z_n $ in $ \mathbb{C}^n $ with $ d_1 c_1 +\cdots + d_n c_n = 0 $ with $ H(z+c)= H(z) $, $ R_1\in\mathbb{C} $. Thus, from \eqref{eq-3.29} we see that
\begin{align*}
\alpha\dfrac{\partial L_1(z)}{\partial z_i} + \beta\dfrac{\partial L_1(z)}{\partial z_j} + \alpha\dfrac{\partial H(z)}{\partial z_i} + \beta\dfrac{\partial H(z)}{\partial z_j} \equiv M_5
\end{align*}
or
\begin{align*}
\alpha\dfrac{\partial H(s)}{\partial z_i} + \beta\dfrac{\partial H(s)} {\partial z_j}\equiv (\alpha d_{i} + \beta d_{j}) H^{\prime}\equiv M_{6},
\end{align*}
where, $ M_5= \frac{2M_1}{M_2}\left(e^{{\eta}/{2}} -1\right) $ and $ M_6 = M_5 - (\alpha a_{1i} + \beta a_{1j}) $. Since $ \alpha d_{i} + \beta d_{j} \neq 0 $, then $ H^{\prime} $ is a constant. Thus, it follows that $ H(s)= A_3 s + A_4 = A_3(d_1 z_1 +\cdots + d_n z_n) + A_4 $, where $ A_3 =\frac{M_4}{\alpha d_{i} + \beta d_{j}}$ and $ A_4 \in\mathbb{C} $. Therefore, we obtain 
\begin{align}\label{eq-3.30}
g(z)= L_1(z) + H(s)+ R_1 = L(z) + R_2 = k_1 z_1 +\cdots + k_nz_n + R_2,
\end{align}
where, $ k_1 =(A_3 d_1 +a_{11}), \ldots, k_n=(A_3 d_n + a_{1n}) $ and $ R_2 = A_4 + R_1 $. In view of \eqref{eq-3.29} and \eqref{eq-3.30}, we obtain
\begin{align*}
	\dfrac{M_2}{2M_1}(\alpha k_i + \beta k_j) + 1 = e^{\frac{k_1 c_1 +\cdots + k_n c_n}{2}}.
\end{align*}
The first equation of \eqref{eq-3.27} can be written as
\begin{align}\label{eq-3.31}
	\alpha\dfrac{\partial f(z)}{\partial z_i} + \beta\dfrac{\partial f(z)}{\partial z_j} = M_{1} e^{\frac{L(z) + R_2}{2}}
\end{align}
Solving the PDE \eqref{eq-3.31}, we obtain
\begin{align}\label{eq-3.32}
	f(z)= \dfrac{2(\omega_2\xi^2 -\omega_1)}{\xi\sqrt{a}(\omega_2-\omega_1)(k_i \alpha + k_j \beta)}e^{\frac{L(z) + R_2}{2}} + \phi_1(z_j -\frac{\beta}{\alpha} z_i).
\end{align}
Moreover, substituting \eqref{eq-3.32} into the second equation of \eqref{eq-3.27} and comparing both sides, we obtain 
\begin{align*}
	\phi_1\left(z_j - \frac{\beta}{\alpha}z_i +c_j -\frac{\beta}{\alpha} c_i\right) = \phi_1\left(z_j - \frac{\beta}{\alpha}z_i\right),
\end{align*}
which implies that $ \phi_1 $ is a finite order periodic entire function with period $ (c_j -\frac{\beta}{\alpha} c_i) $.\\

\noindent{\bf Case B:} If $ e^{h_2(z+c) -h_1(z+c)} $ is not a constant, then obviously, $ H_{21}(z)\equiv 0 $ and $ H_{22}(z)\equiv 0 $ cannot hold at the simultaneously. Otherwise, from \eqref{eq-3.26} we see that $ K_{2} e^{h_2(z+c) -h_1(z+c)}\equiv 1 $, which is a contradiction. \vspace{1.2mm}

If $ H_{21}(z)\equiv 0 $ and $ H_{22}(z)\not\equiv 0 $, in view of \eqref{eq-3.26} we obtain
\begin{align}\label{eq-3.33}
	- H_{22}(z) e^{h_2(z) -h_1(z+c)} + K_{2} e^{h_2(z+c) -h_1(z+c)} \equiv 1,
\end{align}
As $ e^{h_2(z+c) -h_1(z+c)} $ is not a constant, it follows that $ e^{h_2(z) -h_1(z+c)} $ is not a constant. Furthermore, $ e^{h_2(z+c) -h_2(z)} $ is not a constant. Otherwise, $ h_2(z+c) -h_2(z)=\eta_1 $, where $ \eta_1\in\mathbb{C} $. Then, from \eqref{eq-3.33} we see that $ (-H_{22}(z) e^{-\eta_1} + K_2) e^{h_2(z+c) -h_1(z+c)} \equiv 1 $, which is a contraction as $ e^{h_2(z+c) -h_1(z+c)} $ is non-constant. Therefore, the equation \eqref{eq-3.33} can be written as 
\begin{align}\label{eq-3.34}
	- H_{22}(z) e^{h_2(z)} + K_{2} e^{h_2(z+c)} - e^{h_1(z+c)} \equiv 0.
\end{align}
In view of Lemma \ref{lem-3.3}, from \eqref{eq-3.34}, we get contradiction.
\vspace{1.2mm}

Similarly, if $ H_{21}(z)\not\equiv 0 $ and $ H_{22}(z)\equiv 0 $, we get a contradiction. Therefore, we obtain that $ H_{21}(z)\not\equiv 0 $ and $ H_{22}(z)\not\equiv 0 $. As $ h_1(z), h_2(z) $ are polynomials and $ K_{2} e^{h_2(z+c) -h_1(z+c)} $ is non-constant, by Lemma \ref{lem-3.4} for \eqref{eq-3.26}, we obtain
\begin{align*}
	H_{21}(z) e^{h_1(z) -h_1(z+c)}\equiv 1 \;\;\mbox{or}\;\; - H_{22}(z) e^{h_2(z) -h_1(z+c)}\equiv 1.
\end{align*}
\noindent{\bf Sub-case B1:} Assume that $ H_{21}(z) e^{h_1(z) -h_1(z+c)}\equiv 1 $. Then, from \eqref{eq-3.26}, it is easy to see that $ \frac{H_{22}(z)}{K_{2}} e^{h_2(z) -h_2(z+c)}\equiv 1 $. Since $ h_1(z), h_2(z) $ are polynomials, it follows that $ h_1(z) -h_1(z+c)=\eta_2 $ and $ h_2(z) -h_2(z+c)=\eta_3 $, where $ \eta_2, \eta_3\in\mathbb{C} $. Thus, it follows that $ h_1(z)=L_1(z) +H_1(s) + R_3 $ and $ h_2(z)=L_2(z) +H_2(s) + R_4 $, where $ L_{l}(z)= a_{l1} z_1 +\cdots + a_{ln} z_n $ and $ H_{l}(s) $ for $ l=1, 2 $ are polynomial in $ s:= d_1 z_1 +\cdots + d_n z_n $ in $ \mathbb{C}^n $ with $ d_1 c_1 +\cdots + d_n c_n = 0 $ with $ H_{l}(z+c)= H_{l}(z) $ for $ l=1, 2 $, and $ R_3, R_4\in\mathbb{C} $. Since $ \alpha d_{i} + \beta d_{j} \neq 0 $, by the similar argument as in Case $ 1 $ in Theorem \ref{th-2.2}, we see that $ H_l(s) $ is a linear polynomial in $ s $. Therefore, it is easy to see that $ L_{l}(z) + H_l(s) \;(l=1, 2) $ composed of one linear function. For convenience, we always refer to $ h_1(z)=L_1(z) + R_3 $ and $ h_2(z)=L_2(z) + R_4 $. Obviously $ L_1(z)\neq L_2(z) $. Otherwise, $ h_2(z+c) -h_1(z+c) $ becomes a constant, which turns out that $ e^{h_2(z+c) -h_1(z+c)} $ is a constant, a contradiction. Substituting $ h_1(z) $ and $ h_2(z) $ into $ H_{21}(z) e^{h_1(z) -h_1(z+c)}\equiv 1 $ and $ \frac{H_{22}(z)}{K_{2}} e^{h_2(z) -h_2(z+c)}\equiv 1 $, we obtain
\begin{align*}
\begin{cases}
	\dfrac{\sqrt{a}}{\omega_2\sqrt{b}}[(\alpha a_{1i} +\beta a_{1j}) + \sqrt{a}\omega_2] e^{-L_1(c)}\equiv 1 ,\vspace{1.2mm}\\ \dfrac{\sqrt{a}}{\omega_1\sqrt{b}}[(\alpha a_{2i} +\beta a_{2j}) + \sqrt{b}\omega_1] e^{-L_2(c)}\equiv 1.
\end{cases}
\end{align*}
Now, the equation \eqref{eq-3.23} can be written as
\begin{align}\label{eq-3.35}
	\alpha\dfrac{\partial f(z)}{\partial z_i} + \beta\dfrac{\partial f(z)}{\partial z_j} = \dfrac{\omega_2 e^{L_1(z) + R_3} - \omega_1 e^{L_2(z) + R_4}}{\sqrt{a}(\omega_2 -\omega_1)},
\end{align}
solving the PDE \eqref{eq-3.35}, we obtain
\begin{align}\label{eq-3.36}
f(z) =\dfrac{\omega_2 e^{L_1(z) + R_3}}{\sqrt{a} (\omega_2 -\omega_1)(\alpha a_{1i} + \beta a_{1j})} - \dfrac{\omega_1 e^{L_2(z) + R_4}}{\sqrt{a} (\omega_2 -\omega_1)(\alpha a_{2i} + \beta a_{2j})} + \phi_2\left(z_j -\frac{\beta}{\alpha} z_i\right).
\end{align}
Furthermore, substituting \eqref{eq-3.36} into the second equation of \eqref{eq-3.24} and comparing both sides we obtain 
\begin{align*}
	\phi_2\left(z_j - \frac{\beta}{\alpha}z_i +c_j -\frac{\beta}{\alpha} c_i\right) = \phi_2\left(z_j - \frac{\beta}{\alpha}z_i\right),
\end{align*}
which shows that $ \phi_2 $ is a finite order periodic entire function with period $ (c_j -\frac{\beta}{\alpha} c_i) $. From \eqref{eq-3.25}, it follows that
\begin{align*}
	g(z)= h_1(z) + h_2(z) = L(z) + R_5,
\end{align*}
where, $ L(z)=L_1(z) + L_2(z) $ and $ R_5 =R_3 + R_4 $.\vspace{1.2mm}

\noindent{\bf Sub-case B2:} Suppose that $ - H_{22}(z) e^{h_2(z) -h_1(z+c)}\equiv 1 $. Then, from \eqref{eq-3.9} we see that $ -\frac{H_{21}(z)}{K_{2}} e^{h_1(z) -h_2(z+c)}\equiv 1 $. Since $ h_1(z), h_2(z) $ are polynomials, it follows that $ h_2(z) -h_1(z+c)=\eta_4 $ and $ h_1(z) -h_2(z+c)=\eta_5 $, where $ \eta_4, \eta_5\in\mathbb{C} $. A simple computation shows that $ h_1(z+2c) -h_1(z)= -\eta_4 -\eta_5 $ and $ h_2(z+2c) -h_2(z)= -\eta_4 -\eta_5 $. Therefore, we conclude that $ h_1(z)= L(z) + H(s) + R_6 $ and $ h_2(z)= L(z) + H(s) + R_7 $, where $ L(z)= a_1 z_1 +\cdots + a_n z_n $ and $ H(s) $ is a polynomial in $ s:= d_1 z_1 +\cdots + d_n z_n $ in $ \mathbb{C}^n $ with $ d_1 c_1 +\cdots + d_n c_n = 0 $ with $ H(z+c)= H(z) $, and $ R_6, R_7\in \mathbb{C} $. Now, we see that $ h_2(z+c) -h_1(z+c)= R_7 - R_6 $, which shows that $ e^{h_2(z+c) -h_1(z+c)} $ is a constant, a contradiction. This completes the proof.
\end{proof}
\noindent\textbf{Acknowledgment:} The authors would like to thank the referee for their helpful suggestions and comments to improve the exposition of the paper.
\vspace{1.6mm}

\noindent\textbf{Compliance of Ethical Standards:}\\

\noindent\textbf{Conflict of interest.} The authors declare that there is no conflict  of interest regarding the publication of this paper.\vspace{1.5mm}

\noindent\textbf{Data availability statement.}  Data sharing is not applicable to this article as no datasets were generated or analyzed during the current study.

\end{document}